\documentclass[twocolumn,letterpaper]{ieeeconf}

\IEEEoverridecommandlockouts                              
\usepackage{graphicx}
\usepackage{amsmath}
\usepackage{amssymb}
\usepackage{amsfonts}
\usepackage{latexsym}
\usepackage{epsfig}
\usepackage{theorem}
\usepackage{fancybox}
\usepackage{shadow}
\usepackage{float}
\usepackage{comment}
\usepackage{color}
\usepackage{enumerate}
\usepackage{psfrag}

\def \RR {{\mathbb R}}

\def \vec{\overrightarrow}

\def \calI {{\cal I}}
\def \calB {{\cal B}}

\newcommand{\tnm}[1]{\textnormal{#1}}
\newcommand{\ivec}{\vec{\imath}}
\newcommand{\jvec}{\vec{\jmath}}

\begin{document}
\title{\bf Control of VTOL Vehicles with Thrust-direction Tilting}

\author{Minh-Duc Hua, Tarek Hamel, Claude Samson
\thanks{Minh-Duc Hua is with ISIR UPMC-CNRS (Institut des Syst\`emes Intelligents et de Robotique), France. E-mail: $hua@isir.upmc.fr$.}
\thanks{Tarek Hamel is with I3S UNSA-CRNS, France. E-mail: $thamel@i3s.unice.fr$.}
\thanks{Claude Samson is with INRIA, France, and currently on leave at I3S UNSA-CRNS, France. E-mails: $claude.samson@inria.fr$, $csamson@i3s.unice.fr$}}

\maketitle

\pagestyle{empty}
\thispagestyle{empty}

\begin{abstract}
An approach to the control of a VTOL vehicle equipped with complementary thrust-direction tilting capabilities that nominally yield full actuation of the vehicle's position and attitude is developed. The particularity and difficulty of the control problem are epitomized by the existence of a maximal thrust-tilting angle which forbids complete and decoupled control of the vehicle's position and attitude in all situations. This problem is here addressed via the formalism of primary and secondary objectives and by extending a solution previously derived in the fixed thrust-direction case. The proposed control design is also illustrated by simulation results involving a quadrotor UAV with all propellers axes pointing in the same monitored tilted direction.
\end{abstract}

\section{Introduction}
Thrust vectoring for an aerial vehicle is the ability to modify the direction of the propulsion thrust with respect to (w.r.t.) a body-fixed frame. This feature can be used either for attitude (orientation) control, when the thrust rotation center is located at some distance of the vehicle's center of mass (CoM) and thrust vectoring yields torque creation, as in the case of rocket nozzle tilting or ducted-fan airflow derivation via the use of rotating surfaces, or for attitude/position control decoupling, when the thrust rotation center is near the CoM and complementary actuation for attitude control is available, as in the case of VTOL aircraft whose fuselage orientation is controlled independently of the vehicle's longitudinal motion. As a matter of fact, thrust vectoring can also be used to achieve a combination of the aforementioned objectives. This multiple usage renders the term {\em thrust-vectoring} somewhat imprecise. For this reason, we propose here to use the term {\em thrust-tilting} when the second possibility, i.e. attitude and longitudinal motion control decoupling, is addressed, as in the present study. In this case, thrust direction tilting involves two actuated degrees of freedom (d.o.f.) which add to the conventional four actuated d.o.f. associated with common aerial vehicles, i.e. thrust intensity plus three torque components necessary for complete attitude control. This yields six independent actuated d.o.f., thus allowing for the complete control of the six dimensional state associated with the position and attitude of the vehicle's body. A similar objective is addressed in \cite{rbgICRA12} where the vehicle under consideration is a quadrotor UAV whose propellers axes rotate two by two about one of the two orthogonal axes of the quadrotor's ``cross'' structure supporting the propellers. This configuration thus involves four additional motors (one for each propeller) but yields, at the vehicle's body level, only two additional independently actuated d.o.f. as a result of inevitable actuation coupling when more than six independent actuators are employed. The authors of this reference advocate --rightly in our opinion-- that beyond novel challenges that thrust tilting represents from a control point of view, endowing UAVs with full six d.o.f. mobility can be a useful feature in many future applications that involve, for instance, observation and manipulation tasks in cluttered environments. Beyond robotic applications, and with a view to the larger context of aeronautical applications, numerous past and current projects of either manned or unmanned aerial vehicles with thrust tilting also testify of the practical importance of this feature and of mastering its monitoring. Concerning this latter issue, the control design proposed in \cite{rbgICRA12} basically relies on exact linearization of the vehicle's motion equations. This control strategy, combined with actuation redundancy, in turn lead to a control calculation based on the use of pseudo-inverse matrices and on solving a complementary optimization problem (energy expenditure minimization, for example, as proposed in the paper). With respect to this work, we here address thrust tilting in the form of a generic problem whose solution potentially applies to a large panel of aerial vehicles with extended flight envelopes. Such a claim of generality imposes to take aerodynamic forces acting on a moving vehicle into account, even though modelling assumptions in this respect are made here for the sake of simplification, and in order to carry out control solutions for which strong stability and convergence properties can be assessed. While a non-redundant thrust tilting actuation --i.e. involving only the two d.o.f. associated with the modification of an axis direction in the three-dimensional Euclidean space-- is considered the proposed (nonlinear) control design is also different from the one in \cite{rbgICRA12}. It is in fact an extension of the one presented in \cite{hhms09}, \cite{hhms13} for reference velocity or position tracking in the case where the thrust direction is fixed, and it is based on a Lyapunov-like approach. Another important difference and original outcome of this study is that the proposed control solution takes vectoring limitations into account explicitly. More precisely, it is assumed that the thrust tilting angle with respect to a ``neutral'' direction, corresponding for instance to the one associated with the fixed direction of a conventional quadrotor UAV, cannot exceed a known threshold. Due to this limitation, independent control of the vehicle's attitude and position is no longer always possible. Indeed, assume for instance that the main goal is to asymptotically stabilize a given horizontal reference velocity for the vehicle's CoM. Then large reference accelerations and/or strong aerodynamic drag associated with large velocities (relatively to the ambient air) may require a thrust direction inclined almost horizontally. If the maximal tilting angle is small, the necessity of inclining the thrust direction towards the horizontal plane clearly forbids the simultaneous stabilization of any reference attitude associated with hovering, for which the thrust direction is vertical (in the absence of wind). These considerations lead naturally to set priorities as for different control objectives. Accordingly, the proposed control methodology involves a {\em primary} objective associated with the reference velocity or position asymptotic stabilization, and a {\em secondary} objective associated with the asymptotic stabilization of either a reference direction for one of the body-base vectors --thus leaving one d.o.f. unused, or a complete reference attitude for a body-fixed frame. Beside provable stability and convergence properties in a large domain of operation, we believe that the conceptual simplicity of the solution, the non-requirement of switching between several control laws, and the ability of monitoring transition phases smoothly constitute valuable complementary assets. Its geometric nature showing through its construction in the framework of affine geometry and its expression mostly coordinates-free also distinguish it from linear and other nonlinear control methods employed in the domains of aeronautics and aerial robotics.

The remainder of the paper is organized as follows. Section \ref{notation} specifies the notation used along the paper. Section \ref{mod} explains modelling simplifications made at the actuation and control levels, and specifies the class of considered aerodynamic forces acting on the vehicle. Section \ref{controldesign} develops the control design in terms of primary and secondary objectives and shows how to complement feedback control laws derived in previous works in the fixed thrust direction case with a strategy capable of monitoring thrust tilting angle saturation efficiently. Section \ref{application} reports simulation results which illustrate and validate the control approach on a tilted-quadrotor UAV. Section \ref{conclusion} briefly summarizes the contributions of the paper.

\section{Notation} \label{notation}

$\bullet$ With $x \in \RR^n$, $x^T$ stands for the transpose of $x$.

$\bullet$  $\vec{x}$ denotes an affine vector associated with the vector space $\RR^3$.

$\bullet$ The scalar product of two vectors $\vec{x}$ and $\vec{y}$ is denoted as  $\vec{x} \cdot \vec{y}$, and their cross product as $\vec{x} \times \vec{y}$.

$\bullet$ $\vec{x}=\{\vec{\imath},\vec{\jmath},\vec{k}\}x$ is used as a short notation for $\vec{x}=\vec{\imath}x_1+\vec{\jmath}x_2+\vec{k}x_3$, with $x=(x_1,x_2,x_3)^T \in \RR^3$ a vector of coordinates of $\vec{x}$ in the basis $\{\vec{\imath},\vec{\jmath},\vec{k}\}$.

$\bullet$ $G$ : vehicle's center of mass (CoM).

$\bullet$ $\calI=\{O;\vec{\imath_0},\vec{\jmath_0},\vec{k_0}\}$~:~inertial frame.

$\bullet$ $\calB=\{G;\vec{\imath},\vec{\jmath},\vec{k}\}$~:~body-fixed frame.

$\bullet$ $\vec{u}$ : unit vector on the thrust axis.

$\bullet$ $\vec{\omega}^u_{\cal I}=\vec{u} \times \frac{d}{dt} \vec{u}_{| \calI}$: angular velocity of $\vec{u}$ with respect to the inertial frame.

$\bullet$ $\vec{\omega}^u_{\cal B}=\vec{u} \times \frac{d}{dt} \vec{u}_{| \calB}$ : angular velocity of $\vec{u}$ with respect to the body-fixed frame.

$\bullet$ $\vec{\omega}=\{\vec{\imath},\vec{\jmath},\vec{k}\} \omega$: angular velocity of the body-fixed frame w.r.t. the inertial frame.

$\bullet$ $u$: vector of coordinates of $\vec{u}$ in the basis of $\calB$, i.e. $\vec{u}=\{\vec{\imath},\vec{\jmath},\vec{k}\} u$.

$\bullet$ $\dot{u}$: vector of coordinates of $\frac{d}{dt} \vec{u}_{| \calB}$ in the basis of $\calB$, i.e. $\frac{d}{dt} \vec{u}_{| \calB}=\{\vec{\imath},\vec{\jmath},\vec{k}\} \dot{u}$.

$\bullet$ $m$ : vehicle's mass.

$\bullet$ $\vec{F}^a$ resultant of aerodynamic forces acting on the vehicle.

$\bullet$ $\vec{T}=-T \vec{u}$ : thrust force, with the minus sign arising from a convention used for VTOL vehicles.

$\bullet$ $\vec{v}$ : CoM's velocity w.r.t. the inertial frame.

$\bullet$ $\vec{a}$ : CoM's acceleration w.r.t. the inertial frame.

$\bullet$ $\vec{g}$ : gravitational acceleration.
\begin{figure}[htb]
\centering
\includegraphics[width=5 cm, height = 8cm]{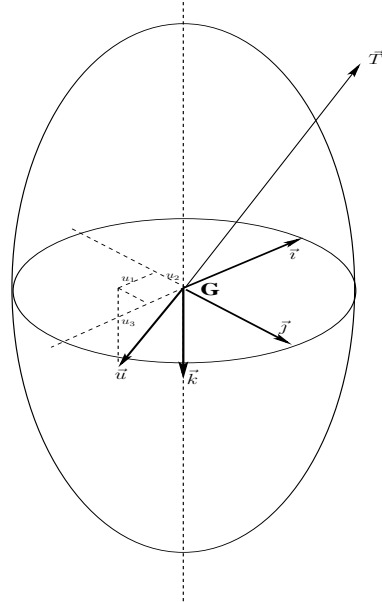}
\vspace{-0.4cm}
\caption{Sketch representation of a VTOL vehicle with thrust tilting.}\label{vtol}
\end{figure}

\section{Modelling simplifications and motion equations}\label{mod}
In this section, a generic control model of a VTOL vehicle evolving in 3D space is proposed. The actuation inputs used to control the vehicle consist of i) a propulsion thrust $\vec{T}$ whose direction with respect to the vehicle's main body can be tilted, ii) a torque $\vec{\Gamma}$, independent of the thrust tilting actuation, used to modify the body's angular velocity $\vec{\omega}$ at will, and iii) torques needed to change the thrust direction. Nominally the resultant thrust force passes through the center of mass and thus exerts zero torque. Tilting the thrust direction may however create a parasitic torque $\vec{\Gamma}_T$. The control torque $\vec{\Gamma}$ then has to pre-compensate for this parasitic torque and also ensure that (almost) any desired angular velocity is physically obtained rapidly (instantaneously, in the ideal case). In practice, this torque is produced in various ways, typically by using propellers (VTOL vehicles \cite{hmlo02}, \cite{rbgICRA12}), or rudders and flaps (aeroplanes \cite{al02}, ducted-fans \cite{nms08}), or control moment gyros (spacecrafts \cite{yt04}), etc. In order to give the reader an idea of how this torque can be calculated to ensure the swift convergence of $\vec{\omega}$ to a desired (time-varying) orientation velocity $\vec{\omega}^\star=\{\vec{\imath},\vec{\jmath},\vec{k}\} {\omega}^\star$, let us consider the classical Euler equation of the orientation dynamics expressed in the coordinates of body-fixed frame:
\begin{equation} \label{euler_equation}
I \dot{\omega}=-\omega \times I\omega + \Gamma_e + \Gamma
\end{equation}
with $I$ the vehicle's angular inertia matrix, $\Gamma_e$ the vector of coordinates in the body-fixed frame ${\cal B}$ of all ``parasitic''
torques (including $\Gamma_T$ and aerodynamic effects) acting on the vehicle. Then, calculating $\Gamma$ as follows
\[\Gamma=-\Gamma_e + I\dot{\omega}^\star + \omega \times I\omega^\star-k_{\omega}I(\omega-\omega^\star)~~;~k_{\omega}>0 \]
and applying this torque to the vehicle yields the closed-loop equation
\[I\frac{d}{dt}(\omega-\omega^\star)=-\omega \times I(\omega-\omega^\star)-k_{\omega}I(\omega-\omega^\star)\]
and hence the exponential stability of $(\omega-\omega^\star)=0$ with a rate of convergence of $\omega$ to $\omega^\star$ given by $k_{\omega}$. Under some extra assumptions upon $\omega^\star$ and $\Gamma_e$, one shows that a simple proportional feedback $\Gamma=-Ik_{\omega}(\omega-\omega^\star)$ with $k_{\omega}>0$ large enough, suffices to render and maintain $|\omega-\omega^\star|$ smaller than any given small threshold.

The above considerations justify the conceptual ``backstepping'' assumption, commonly made, which consists in using the angular velocity $\vec{\omega}$, rather than the torque $\vec{\Gamma}$ which produces this velocity, as a control input. In the same order of idea, we will assume that the thrust intensity $T$, and the thrust-direction tilting angular velocity $\vec{\omega}_{\cal B}^{u}$, are the other control inputs at our disposal. The equations characterizing the system's dynamics are then
\begin{align}
m \vec{a} &= m \vec{g} + \vec{F}^a - T \vec{u} \label{acc0} \\
\frac{d}{dt} \{\vec{\imath},\vec{\jmath},\vec{k}\}_{| \calI}&= \vec{\omega} \times \{\vec{\imath},\vec{\jmath},\vec{k}\} \label{omegab}\\
\frac{d}{dt} \vec{u}_{| \calB}&=\vec{\omega}^u_{\cal B} \times \vec{u} \label{omegau}
\end{align}
with $\vec{F}_a$ representing all forces, other than thrust and gravitation, applied to the vehicle's body. This vector is typically dominated by lift and drag aerodynamic forces whose intensities depend on the vehicle's longitudinal velocity (relatively to the ambient air) and, of particular importance for control purposes, on the vehicle's attitude when lift is not negligible. For the sake of simplification, and because the focus of the present paper is not to discuss control design aspects specifically related to the aerodynamic forces acting on the vehicle, we will here assume that $\vec{F}^a$ can be decomposed into the sum of two components as follows
\begin{equation}\label{Fa}
\vec{F}^a = \vec{F}^a_1  + \vec{F^a_2}
\end{equation}
with $\vec{F^a_2}=F^a_2 \vec{u}$ ($F^a_2 \in \RR$) and $\vec{F}^a_1$ not depending on the vehicle's orientation. Note that this assumption is not as restrictive as it may first seem. For instance, if lift is negligible (the case of a body with spherical shape), then $\vec{F}^a$ is essentially reduced to the drag force component $\vec{F}^a_1=-c_D|\vec{v}^a|\vec{v}^a$, where $c_D$ is the drag characteristic coefficient, $\vec{v}^a :=\vec{v}-\vec{v}^w$ is the apparent velocity (i.e. the vehicle's velocity w.r.t. the ambient air), and $\vec{v}^w$ is the wind velocity w.r.t. an inertial frame. In the case of a quadrotor UAV with tilted thrust-direction, alike the one considered in the simulation Section \ref{application}, the expression of this force is more complex. It involves a body-drag force $\vec{F}^a_D$ and an induced-drag force $\vec{F}^a_I$ generated by the airflow circulation around the rotors blades. Expressions of these forces are \cite{mkc12}:
\[
\vec{F}^a_D = -c_D |\vec{v}^a| \vec{v}^a, \;
\vec{F}^a_I = -c_I (\vec{ v}^a - (\vec v^a \cdot \vec u) \vec{ u})
\]
with $c_I$ an aerodynamic coefficient. Summing up these two forces yields
\begin{align*}
\vec{F}^a &= \vec{F}^a_D  + \vec{F}^a_I \\
          &=\left(-c_D |\vec{v}^a| \vec{v}^a -c_I \vec{v}^a\right ) +  \left ( c_I(\vec{v}^a \cdot \vec{ u} ) \vec{u}\right)\\
          &=\vec{F}^a_1  + \vec{F^a_2}
\end{align*}
with $\vec {F}^a_1:= -c_D |\vec{v}^a| \vec{v}^a -c_I \vec{v}^a$ independent of the vehicle's orientation, and $\vec F^a_2:=  c_I(\vec{v}^a \cdot \vec{ u} ) \vec{u}$ aligned with the thrust direction $\vec u$, so that the above assumption is again satisfied. This assumption also holds in the first approximation, and in a large domain of operation, for all bisymmetric vehicles with fixed thrust direction, as shown in \cite{pucci12-rap}.

In view of \eqref{Fa}, the longitudinal dynamics equation \eqref{acc0} can be rewritten as
\begin{equation} \label{acc}
m \vec{a} = m \vec{g} +\vec{F}^a_1 - \bar{T} \vec{u}
\end{equation}
with $\bar{T}=T-F^a_2$ and an external force $\vec{F}^a_1$ which does not depend on the vehicle's orientation. We also assume that ``reasonably'' good on-line knowledge of $\vec{F}^a$ and $\vec{{T}}$ is available for the control calculations. Control robustness should compensate for moderate modelling errors. However, to simplify the presentation of the control design, it is convenient to assume that these terms are known exactly.

\section{Control design} \label{controldesign}
In the previous section we have justified the use of $\bar{T}$, $\vec{\omega}$, and $\vec{\omega}^u_{\cal B}$ as control inputs. The present section aims at working out feedback control expressions for these inputs.

\subsection{Primary objective realization}
From the definitions of the angular velocities of $\vec{\omega}$, $\vec{\omega}^u_{\cal I}$, and $\vec{\omega}^u_{\cal B}$, one has
\begin{equation} \label{basics}
\vec{\omega}^u_{\cal I}=\vec{\omega}^u_{\cal B}+\vec{\omega}-\vec{u}(\vec{u} \cdot \vec{\omega})
\end{equation}

We can assume that $\bar{T}$ and $\vec{\omega}^u_{\cal I}$ are determined in order to achieve a primary objective related to the vehicle's longitudinal motion. For instance, consider the case where velocity reference asymptotic stabilization is the objective. Denoting the velocity error as $\vec{\tilde{v}}=\vec{v}-\vec{v}_r$, one obtains the following error equation
\begin{equation}
m\frac{d}{dt} \vec{\tilde{v}} |_{\cal I}=m(\vec{a}-\vec{a}_r) = - \bar{T} \vec{u} + \vec{F} \label{syst-hhms1}
\end{equation}
with
\begin{equation}
\vec{F} := m\vec{g} + \vec{F}^a_1 - m \vec{a}_r
\end{equation}
where $\vec{a}_r=\frac{d}{dt} \vec{v}_{r~| \calI}$ represents the acceleration of the reference trajectory w.r.t. the inertial frame.

Let ${T}_r$ and $\vec{u_r}$ denote the (positively signed) intensity and direction vector of $\vec{F}$ respectively, i.e.
\begin{equation} \label{Trur}
({T}_r, \vec{u}_r)  \equiv \left(|\vec{F}|,\frac{\vec{F}}{|\vec{F}|} \right)
\end{equation}
As long as $|\vec{F}|\neq 0$, defining these terms poses no problem. We will assume from now on that this condition, which is automatically satisfied in the case of hovering thanks to the gravitational acceleration, is satisfied. As explained in \cite{hhms09}, the satisfaction of this condition is also necessary to the controllability of the linearized error system and the existence of conventional, either linear or nonlinear, control solutions that do not rely on the execution of manoeuvres performed in 3D-space. In view of \eqref{syst-hhms1}, one must have $\bar{T} \vec{u}=\vec{F}$ at the equilibrium $\vec{\tilde{v}}=0$. Therefore $\bar{T}$ must converge to $T_r$, whereas the thrust direction $\vec{u}$ must converge to $\vec{u}_r$.

The above considered problem of asymptotic stabilization of a reference velocity has already received solutions (see \cite{hhms09} or \cite{hhms13}) which, modulo minor modifications and transposition details, yield the following control expressions
\begin{align}
\bar{T} &= {T}_r \vec{u} \cdot \vec u_r + k_1  m \vec{u} \cdot \vec{\tilde{v}} \label{Tbar} \\
\vec{\omega}^u_{\cal I} & = \frac{k_2 m}{{T}_r} \vec{u} \times \vec{\tilde{v}}+
\frac{(k_3+\bar{k}_3)}{(1+\vec{u} \cdot \vec{u}_r)^2} \vec{u} \times \vec{u}_r \notag \\
 ~ & \quad - \vec{u} \times \left (\vec{u} \times \left (\vec{u}_r \times \frac{d}{dt}\vec{u}_{r| \cal I} \right )\right)
\label{omega}
\end{align}
with  $k_1$, $k_2$ and $k_3$ positive gains, and $\bar{k}_3=2\dot{T}_r(1+ \vec u \cdot \vec u_r)/T_r$. This controller is derived by considering the following Lyapunov function candidate:
\[
{\cal L}=\frac{1}{k_2}\frac{T_r^2}{m^2} \left (1-\vec{u} \cdot \vec{u}_r\right ) + \frac{1}{2}|\vec{\tilde v}|^2
\]
whose time-derivative along any solution to the controlled system is:
\[
\dot{\cal L}= - \frac{k_3}{k_2} \frac{T_r^2}{m^2} \frac{\left | \vec{u} \times \vec{u}_r \right |^2}{\left (1+ \vec{u} \cdot \vec{u}_r\right)^2} - k_1 (\vec{u} \cdot \vec{\tilde{v}})^2~(\leq 0)
\]
Provided that $|\vec{F}|\neq 0$, that the reference velocity and acceleration are bounded, and that $\vec{u}+\vec{u}_r \neq 0$ initially, then $(\vec{\tilde{v}}=\vec{0},\vec{u}-\vec{u}_r=\vec{0})$ is asymptotically stable (see \cite{hhms09} for the proof of this result).

If, instead of velocity reference stabilization, the primary objective is the tracking of a reference position trajectory, the same control can be used modulo minor modifications involving (bounded) integrals of velocity and position errors. The Lyapunov function candidate used for convergence and stability analysis is modified accordingly \cite{hhms09}, \cite{hhms13}.
For the sake of precision, let $\vec{x}:= \vec{OG}$ and $\vec{x_r}$ respectively denote the vehicle's position and reference position w.r.t an inertial frame, and let $\vec{\tilde{x}}:=\vec{x}-\vec{x}_r$ denote the position error vector.
The above mentioned modifications consist, for instance, in replacing the definition \eqref{Trur} of $T_r$ and $u_r$ by
\begin{equation} \label{Trurxi}
({T}_r, \vec{u}_r)  \equiv \left(|\vec{F}_{\xi}|,\frac{\vec{F_{\xi}}}{|\vec{F_{\xi}}|} \right)
\end{equation}
where:

$\bullet$ $\vec F_{\xi} := m\vec g + \vec F^a_1 - m \vec a_r + m k_I \frac{d^2}{dt^2}\vec z_{| \cal I} + m \vec\sigma(\vec{\tilde \xi})$

$\bullet$ $\vec{{\tilde \xi}} := \vec{\tilde x} + k_I \vec z ~~~;~k_I>0$

$\bullet$ $\vec z$ is a bounded integral of the position error $\vec{\tilde{x}}$ defined as the solution to the following equation
\[
\begin{split}
\frac{d^2}{dt^2}\vec z_{| \cal I} \!=\!&-k_{\dot{z}} \frac{d}{dt}\vec z_{| \cal I}
\\&+ \tnm{sat}^{\frac{\ddot z_{max}}{2}}(k_{z}(-\vec z + \vec{\tnm{sat}}^{\Delta_z}(\vec z + {\vec{\tilde x}}/{k_{z}})))
\end{split}
\]
with zero initial conditions, i.e. $\vec z(0)  = \frac{d}{dt}\vec z_{| \cal I}(0) = \vec 0$, positive numbers $k_{z}$, $k_{\dot{z}}$, $\ddot z_{max}$ and $\Delta_z$, and $\vec{\tnm{sat}}^\Delta$ denoting the classical saturation function defined by $\vec{\tnm{sat}}^\Delta(\vec x) = \min(1,\Delta/|\vec x|) \vec x$.

$\bullet$ $\vec\sigma$ is a bounded function satisfying some properties (detailed in  \cite{hhms09} and \cite{hhms13}) and an example of which is the function defined by $\vec\sigma(\vec{y}) := \beta (\beta^2 |\vec{y}|^2/\eta^2+ 1)^{-\frac{1}{2}} \vec{y}$,
with $\beta$ and $\eta$ denoting positive numbers.

The last modification concerns the term $\vec{\tilde{v}}$ which, in the control expressions \eqref{Tbar} and \eqref{omega}, has to be replaced by $\vec{\tilde{v}}_\xi := \vec {\tilde{v}} + k_I \frac{d}{dt}\vec{z}_{|\calI}$.

\subsection{Secondary objective realization}
The thrust intensity $T$ (or equivalently, $\bar{T}$) has been determined previously (the relation \eqref{Tbar}). It thus remains to determine $\vec{\omega}$ and $\vec{\omega}^u_{\cal B}$ that satisfy \eqref{basics}, as imposed by the realization of the primary objective, and also allow for the realization of a secondary objective.

Let ${\vec{\omega}}^\star$ denote the angular velocity control that would be used for the secondary control objective if the thrust tilting angle was not limited. For instance, this objective can be the asymptotic stabilization of the body-fixed frame vector $\vec{k}$ at a desired time-invariant unit vector $\vec{\eta}$. Then, a possible control whose expression is reminiscent of the one of $\vec{\omega}_{\cal I}^u$ in \eqref{omega} is
\begin{equation} \label{omegabstar1}
{\vec{\omega}}^\star= \frac{k_4}{(1+\vec{k}.\vec{\eta})^2} \vec{k} \times \vec{\eta} + \lambda \vec{k}~;\quad k_4>0,~\lambda \in \RR
\end{equation}
The indetermination of $\lambda$ corresponds to the unused d.o.f associated with rotations about the $\vec{k}$ direction. Now, if the secondary objective is to asymptotically stabilize a desired reference attitude for the vehicle's body, then a possible control is
\begin{equation} \label{omegabstar2}
{\vec{\omega}}^\star=-k_4 \tan (\theta/2) \vec{\nu}~; \quad k_4>0
\end{equation}
with $\theta \vec{\nu}$ the rotation vector associated with the attitude error between the body-fixed frame and the desired orientation.

Let us define
\begin{equation} \label{omegaustar}
{\vec{\omega}^u_{\cal B}}^\star:=\vec{\omega}^u_{\cal I}-({\vec{\omega}}^\star-\vec{u}(\vec{u} \cdot {\vec{\omega}}^\star))
\end{equation}
Since ${\vec{\omega}^u_{\cal B}}^\star$ satisfies \eqref{basics} when $\vec{\omega}$ is equal to the unconstrained solution ${\vec{\omega}}^\star$, it is the thrust tilting angular velocity that would be used in the unconstrained case to achieve both primary and secondary objectives. However, due to the thrust tilting angle limitation, modified expressions for $\vec{\omega}$ and $\vec{\omega}^u_{\cal B}$ have to be worked out.
To this aim, let us first specify the maximal value of the thrust tilting angle and set
\begin{equation} \label{delta}
\delta = \max\left(\sqrt{u_1^2+u_2^2}\right)~~~(<1)
\end{equation}
so that the maximal tilting angle is equal to $\arcsin(\delta)$ ($\in (0,\frac{\pi}{2})$.

Define
\begin{equation} \label{dvecustar}
{\dot{\vec{u}}}^\star:={\vec{\omega}^u_{\cal B}}^\star \times \vec{u}
\end{equation}
and denote ${\dot{{u}}}^{\star}$ as the vector of coordinates of ${\dot{\vec{u}}}^{\star}$ in the basis of the body-fixed frame $\calB$, i.e. ${\dot{\vec{u}}}^{\star}=\{\vec{\imath},\vec{\jmath},\vec{k}\}{\dot{{u}}}^{\star}$. Using $x_{1,2}$ to denote the vector of first two components of $x \in \RR^3$, we propose to modify the tilting angle according to the following control law:
\begin{equation} \label{du12}
\left\{
\begin{array}{lll}
\dot{u}_{1,2} & = &-k_u {u}_{1,2} +k_u \tnm{sat}^{\delta}({u}_{1,2}+ {\dot{u}_{1,2}^\star}/{k_u}) \\
\dot{u}_3 & = & -{u}_{1,2}^T\dot{u}_{1,2}/u_3
\end{array}
\right.
\end{equation}
with $k_u$ a positive number (not necessarily constant), and $\tnm{sat}^{\delta}(.)$ the classical saturation function defined by $\tnm{sat}^{\delta}(x)=\min(1,{\delta}/{|x|})x$, with $x \in \RR^2$ in the present case. Then, $\vec{\omega}^u_{\cal B}$ is given by
\[
\vec{\omega}^u_{\cal B} = \vec u \times \frac{d}{dt} \vec u_{| \cal B}
\]
with $\frac{d}{dt} \vec u_{| \cal B} = \{\ivec ,\jvec,\vec k\} \dot u$. One easily verifies from \eqref{du12} that $\dot{u}=\dot{u}^\star$ as long as $|{u}_{1,2}+{\dot{u}_{1,2}^\star}/{k_u}| \leq \delta$. This indicates that, for ${u}_{1,2}$ to be dominant in the left-hand side term, $k_u$ should be chosen large enough. In view of \eqref{du12}, one also has the relations
\[
\begin{array}{lll}
\frac{1}{2} \frac{d}{dt} |{u}_{1,2}|^2 & = & -k_u |{u}_{1,2}|^2 +k_u {u}_{1,2}^T \tnm{sat}^{\delta}({u}_{1,2}+ {\dot{u}_{1,2}^\star}/{k_u}) \\
~ & \leq & -k_u |{u}_{1,2}|^2+k_u \delta |{u}_{1,2}|
\end{array}
\]
from which one deduces that $|{u}_{1,2}|$ remains smaller or equal to $\delta$, provided that the initial value of $|{u}_{1,2}|$ is itself chosen smaller or equal to $\delta$. Moreover, using the fact that, by definition of ${\dot{\vec{u}}}^{\star}$ (see \eqref{dvecustar}),
$\dot{u}_3^\star=-{u}_{1,2}^T \dot{u}_{1,2}^\star/u_3$, it comes that $\dot{u}=\dot{u}^\star$ whenever $|{u}_{1,2}+\dot{u}_{1,2}^\star/k_u| \leq \delta$.
Therefore, this control law i) respects the imposed limitation on the thrust tilting angle, and ii) allows for the realization of the secondary objective as long as $|{u}_{1,2}+\dot{u}_{1,2}^\star/k_u| \leq \delta$.

It then only remains to determine ${\vec{\omega}}$. An obvious choice is
\begin{equation} \label{omegabfinal}
\vec{\omega} = \vec{\omega}_{\cal I}^u -\vec{\omega}^u_{\cal B}+\vec{u}(\vec{u} \cdot {\vec{\omega}}^\star)
\end{equation}
Indeed, one deduces from this relation that $\vec{u} \cdot {\vec{\omega}}=\vec{u} \cdot {\vec{\omega}}^\star$ and, subsequently, that $\vec{\omega}^u_{\cal B}={\vec{\omega}^u_{\cal B}}^\star$ and $\vec{\omega}={\vec{\omega}}^\star$ when $|{u}_{1,2}+ {\dot{u}_{1,2}^\star}/{k_u}| \leq \delta$. Moreover, the equality \eqref{basics} is always satisfied with this choice, in accordance with the priority given to the realization of the primary objective.

\section{Application to a thrust-tilted quadrotor} \label{application}

\begin{figure}[htb]
\centering
\includegraphics[width=0.9\linewidth]{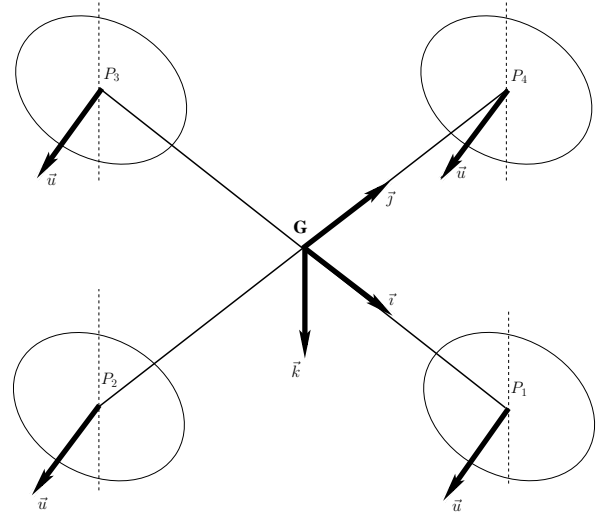}
\caption{Sketch representation of a quadrotor with thrust tilting.}\label{quadrotor}
\end{figure}
The above control solution has been tested on the model of a quadrotor UAV sketched on Fig. \ref{quadrotor}. We assume that the four rotor axes can be tilted and point in the same direction corresponding to the overall thrust direction $\vec{u}$. Let $P_{i}$ ($i=1,\ldots,4$) be the pivoting points of the four rotors, with their positions in the body frame defined by
\[
\begin{split}
\vec {GP}_1 &= h \vec k + d \ivec, \quad \vec {GP}_2 = h \vec k - d \jvec \\
\vec {GP}_3 &= h \vec k - d \ivec, \quad \vec {GP}_4 = h \vec k + d \jvec
\end{split}
\]
with  $h \in \RR$ and $d>0$. On Fig. \ref{quadrotor}, $h=0$.
Let $\varpi_i$ ($i=1,\ldots,4$) denote the angular velocities of the four rotors. According to \cite{hmlo02}, the i-th rotor, generates a thrust force $\vec T_i = \mu \varpi_i^2 \vec u$ and a drag torque $\vec Q_i = \lambda_i \kappa \varpi_i^2 \vec u$, with $\mu$ and $\kappa$ two aerodynamic coefficients and $\lambda_i = 1$ (resp. $\lambda_i = -1$) if $i$ is odd (resp. even).
The thrust $T$ and the torque vector $\vec \Gamma = (\ivec ,\jvec ,\vec k)\Gamma$ generated by the four rotors are thus given by
\[
\begin{split}
T& = \sum_i T_i= \mu  \left(\sum_i\varpi_i^2\right) \vec u \\
\vec \Gamma &= \sum_i (\vec {GP}_i  \times \vec T_i + \vec Q_i)\\&= d \mu (\varpi_1^2 - \varpi_3^2) (\ivec \times \vec u) -d \mu (\varpi_2^2 - \varpi_4^2) (\jvec \times \vec u)\\
&\quad \quad \quad \quad \quad \quad +h \mu \sum_i \varpi_i^2 (\vec k \times \vec u) + \sum_i \lambda_i \kappa \varpi_i^2 \vec u
\end{split}
\]

One deduces the following relation between the vector of rotors angular velocities $\varpi_i$ on the one hand, and the vector composed of the thrust intensity and the control torque components on the other hand
\[
\begin{bmatrix} T \\ \Gamma \end{bmatrix} = A_{mot} \begin{bmatrix} \varpi_1^2 \\ \varpi_2^2 \\ \varpi_3^2\\ \varpi_4^2 \end{bmatrix}
\]
with the allocation matrix $A_{mot}  = [a_{c1} \quad a_{c2} \quad a_{c3} \quad a_{c4}]$,
\[
\begin{split}
&a_{c1}\!:=\! \begin{bmatrix}\mu \\  -h\mu u_2 \!+\! \kappa u_1 \\ h\mu u_1 \!-\! d\mu u_3\!+\!\kappa u_2 \\  d\mu u_2\! +\! \kappa u_3 \end{bmatrix}\!\!,
a_{c2}\!:=\! \begin{bmatrix}\mu \\  -\!h\mu u_2\!-\!d\mu u_3 -\!\kappa u_1 \\ h\mu u_1 \!-\!\kappa u_2 \\ \!d\mu u_1\! -\!\kappa u_3 \end{bmatrix}
\end{split}
\]
\[
\begin{split}\vspace{0.5cm}
&a_{c3}\!:=\! \begin{bmatrix}\mu \\ -\!h\mu u_2\!+\!\kappa u_1 \\ h\mu u_1\!+\!d\mu u_3\!+\! \kappa u_2   \\  -d\mu u_2\!+\!\kappa u_3  \end{bmatrix}\!\!,
a_{c4}\!:=\! \begin{bmatrix}\mu \\ -\!h\mu u_2\! +\!d\mu u_3\!-\!\kappa u_1   \\h\mu u_1 \!-\!\kappa u_2   \\  \!-d\mu u_1 \!-\!\kappa u_3  \end{bmatrix}
\end{split}
\]
Since $A_{mot}$ is invertible ($\text{det}(A_{mot} ) = 8 \kappa d^2 \mu^3 u_3 >0$), $T$ and $\Gamma$ can be given any desired values --modulo the constraint of positivity of the rotors angular velocities and the limited range of velocities imposed by power limitations of the rotors--, and can thus be used as independent control variables. The direct application of the proposed control strategy relies on this actuation property.

\subsection{Simulation results}
Specifications of the simulated vehicle are given in Tab. \ref{tab:specs}.

\begin{table}[h]
\caption{Specifications of the simulated quadrotor} \label{tab:specs} \vspace{-0.4cm}
\begin{center}
\begin{tabular}{|c||c|}
\hline
Specification & Numerical Value\\
\hline
Mass $m$ [$kg$] & 1.5\\
Moment of Inertia $I$ [$kg$ $m^2$] & diag(0.028,0.028, 0.06)\\
Level arm values $[h,d]$ $[m]$ & [0.05,0.2] \\
Thrust angle limitation  $[rad]$ & $\pi/6$ \\
Body drag coefficient $c_{D}$ [$kg$ $m^{-1}$]& 0.0092\\
Induced drag coefficient $c_{I}$ [$kg$ $s^{-1}$]& 0.025\\
\hline
\end{tabular}
\end{center}
\end{table}
Concerning the calculation of the torque $\Gamma$ in charge of producing the desired body angular velocity defined by \eqref{omegabfinal}, we have used $\Gamma = -k_\omega I(\omega - \omega^{\star}) + \omega \times I \omega^{\star}$, with $\omega$ denoting the vehicle's angular velocity obtained by integration of the Euler equation \eqref{euler_equation}, $\omega^{\star}$ the vector of coordinates, expressed in the body-fixed frame $\calB$, of the reference angular velocity defined by \eqref{omegabfinal}, and $k_\omega$ a positive gain. Due to the parasitic torque induced by the chosen non-zero value of $h$ (one of the parameters characterizing the position of the propellers), the non pre-compensation of both this torque and the angular acceleration $\dot{\omega}$ in the expression of $\Gamma$, there remains a residual error between $\omega$ and $\omega^{\star}$ in the general case. This is a deliberate choice to test the robustness of the proposed control design w.r.t. (inevitable) modelling errors, and also the reason of the residual position tracking errors that can be observed from the simulation results reported further on.

The primary objective considered for these simulations is the position tracking of an eight-shaped Lissajous trajectory defined by
\begin{equation}\label{reference}
\vec x_{r} = 5 \sin(a_r t) \vec \imath_0 +   5 \sin(2 a_r t) \vec \jmath_0
\end{equation}
This is a closed curve. By playing on the choice of the parameter $a_r$, one modifies the time period of a complete run, as well as the associated reference velocity and acceleration. The larger $a_r$, the shorter the time period, and the larger the reference velocity and acceleration in the average. Two values of $a_r$ ($2 \pi/15$ and $\pi /5$) are considered. The first one corresponds to a ``slow'' run (Simulation 1) that involves non-saturated thrust-direction tilt angles, whereas the second one corresponds to a ``fast'' run (Simulation 2) along portions of which the tilt angle reaches its maximal value.

The chosen secondary objective is the stabilization of the vehicle's attitude about the identity matrix. In particular, the realization of this objective implies that the vehicle's plane containing the rotors pivoting points remains horizontal all the time. The associated angular velocity control $\vec{\omega}^\star$ is given by \eqref{omegabstar2}.

\begin{figure}[htb]\centering
\psfrag{x(m)}{\scriptsize $x_1\, (m)$}%
\psfrag{y(m)}{\scriptsize $x_2\, (m)$}%
\psfrag{reference trajectory}{\scriptsize reference trajectory}%
\psfrag{vehicle's trajectory}{\scriptsize vehicle's trajectory}%
\includegraphics[width=8.5 cm,  height=5.5cm]{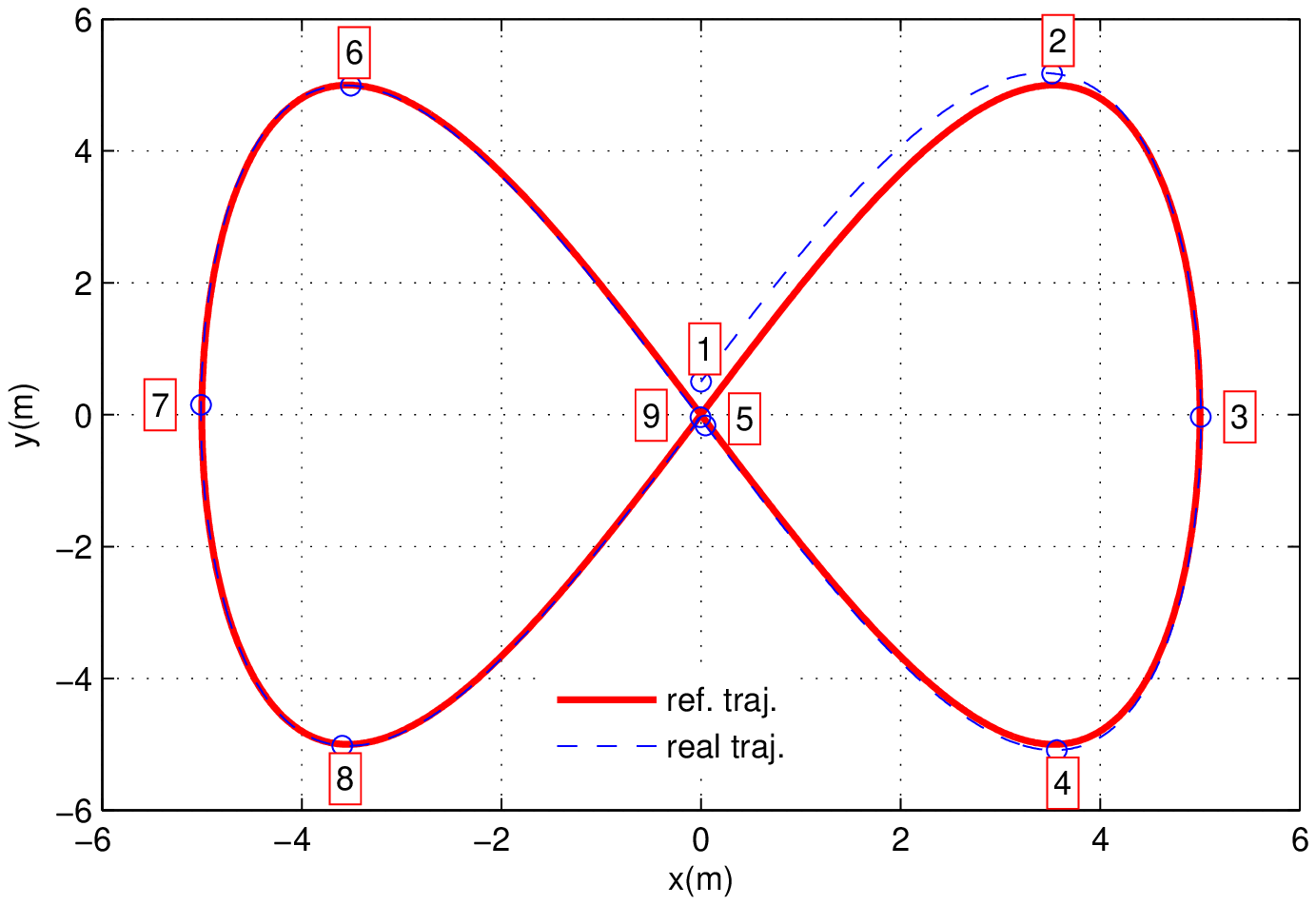}
\vspace*{-0.3cm}
\caption{Reference and vehicle trajectories projected on the horizontal plane (Simulation 1).} \label{fig:simu1a}
\vspace*{0.2cm}
\psfrag{t (s)}{\scriptsize $t\, (s)$}%
\psfrag{angle (deg)}{\scriptsize $Inclination$ $angle$ $(deg)$}%
\psfrag{thrust orientation}{\scriptsize thrust tilt angle}%
\psfrag{vehicle orientation}{\scriptsize vehicle's inclination angle}%
\includegraphics[width=8.5 cm, height=5.5cm]{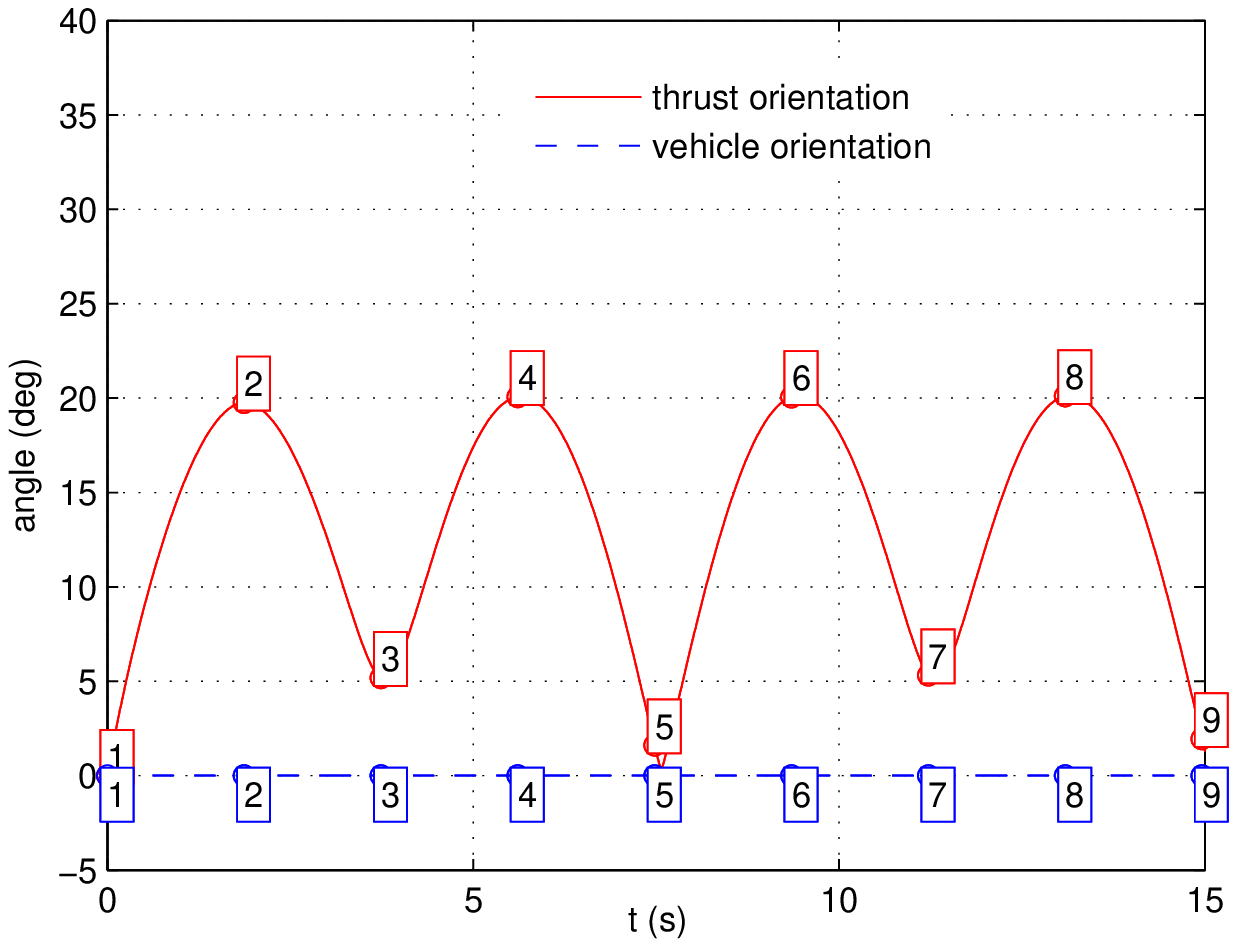}
\vspace*{-0.3cm}
\caption{Thrust and vehicle inclination angles vs. time (Simulation 1).} \label{fig:simu1b} \vspace*{0.2cm}
\psfrag{t (s)}{\scriptsize $t\, (s)$}%
\psfrag{position error (m)}{\scriptsize $Position$ $errors$ $\tilde{x}$ $(m)$}%
\psfrag{error x}{\scriptsize $\tilde{x}_1(m)$}%
\psfrag{error y}{\scriptsize $\tilde{x}_2(m)$}%
\psfrag{error z}{\scriptsize $\tilde{x}_3(m)$}%
\psfrag{vehicle orientation}{\scriptsize vehicle's inclination angle}%
\includegraphics[width=8.5 cm, height=5.5cm]{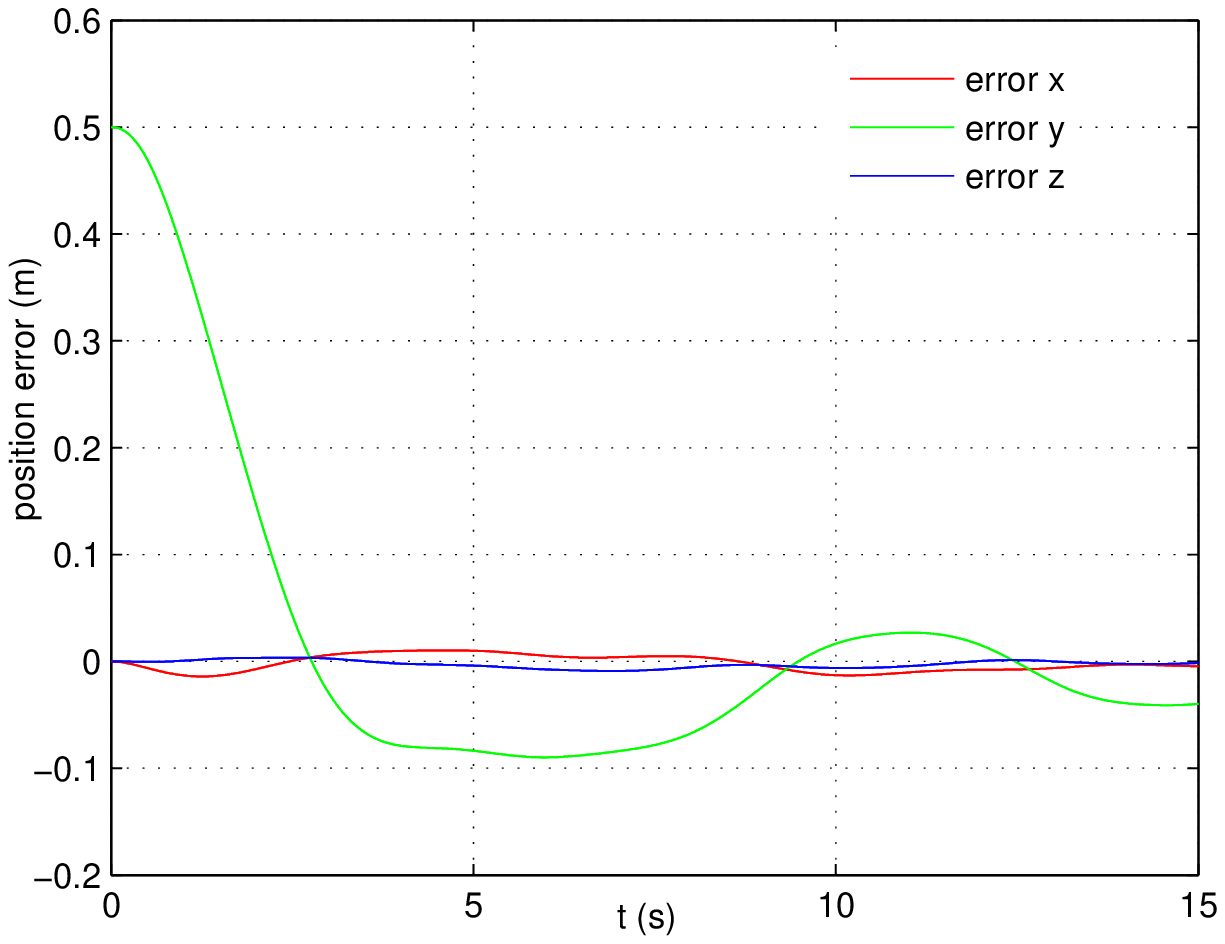}
\vspace*{-0.3cm}
\caption{Position tracking errors vs. time (Simulation 1).} \label{fig:simu1c} \vspace*{0.cm}
\end{figure}

\begin{figure} \centering
\psfrag{t (s)}{\scriptsize $t\, (s)$}%
\psfrag{Thrust (N)}{\scriptsize $Thrust$ $T(N)$}%
\includegraphics[width=8.5 cm, height=5.5cm]{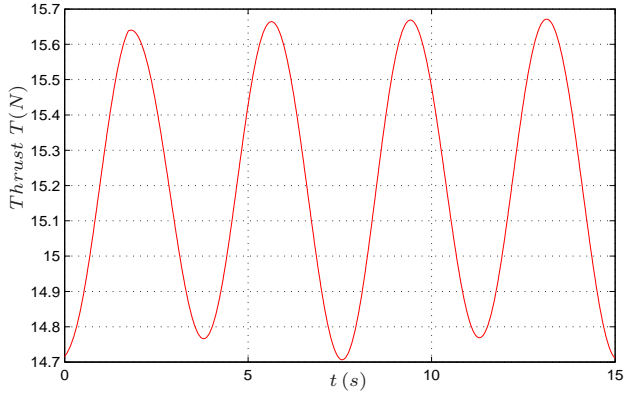}
\vspace*{-0.3cm}
\caption{Thrust intensity vs. time (Simulation 1).} \label{fig:simu1d} 
\end{figure}

The control gains and other parameters involved in the calculation of the control inputs are chosen as follows:

\begin{itemize}
\item $k_1 = 1.2$, $k_2 = 0.34$, $k_3 = 12.8$, $k_I = 1$,
\item $\beta  = 0.36$, $\eta  = 6$,
\item $k_{\dot{z}} = 4$, $k_{z}=4$, $\Delta_{z} = 1$, $\ddot z_{max} =  0.5$,
\item $k_4 = 10$, $k_u = 16$,
\item $k_\omega = 20$.
\end{itemize}

Initial conditions for the vehicle's configuration are as follows:
\[
\left\{
\begin{split}
&\vec x(0)= 0.5 \vec \jmath_0,\quad \vec v(0)= 5a_r \vec \imath_0 + 10 a_r\vec \jmath_0, \\
& \{\ivec (0), \jvec (0), \vec k(0)\}= \{\ivec_0, \jvec_0, \vec k_0\},\\
&  \vec u(0)= \vec k(0), \quad\vec \omega(0) =\vec \omega^u_{\cal B}(0) = \vec 0.
\end{split} \right.
\]

Simulation results are reported next.

\begin{figure}[htpb!]\centering%
\psfrag{x(m)}{\scriptsize $x_1\, (m)$}%
\psfrag{y(m)}{\scriptsize $x_2\, (m)$}%
\psfrag{reference trajectory}{\scriptsize reference trajectory}%
\psfrag{vehicle's trajectory}{\scriptsize vehicle's trajectory}%
\includegraphics[width=8.5 cm,  height=5.5cm]{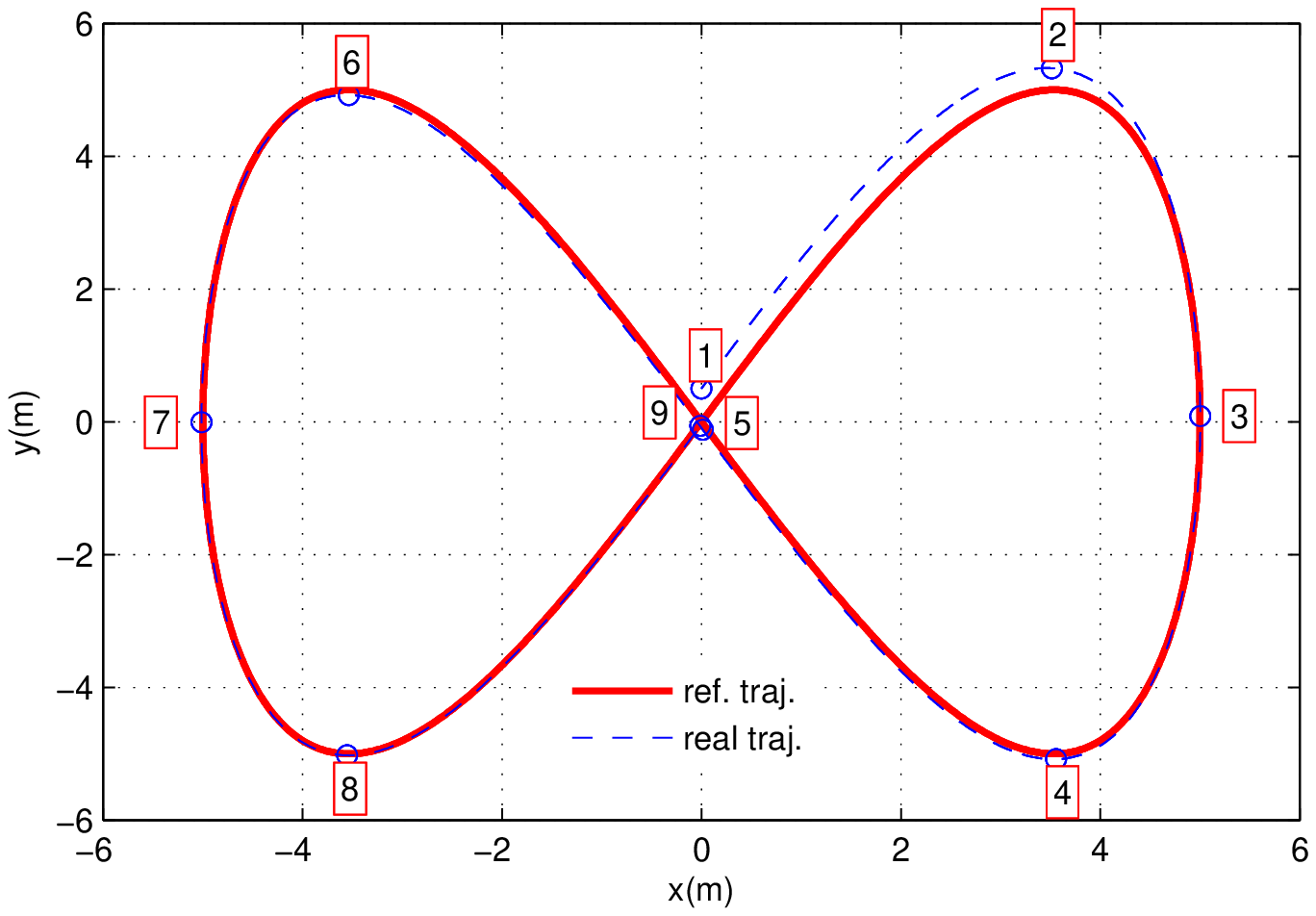}
\vspace*{-0.3cm}
\caption{Reference and vehicle trajectories projected on the horizontal plane (Simulation 2).} \label{fig:simu2a}
\vspace*{0.2cm}
\psfrag{t (s)}{\scriptsize $t\, (s)$}%
\psfrag{angle (deg)}{\scriptsize $Inclination$ $angle$ $(deg)$}%
\psfrag{thrust orientation}{\scriptsize thrust tilt angle}%
\psfrag{vehicle orientation}{\scriptsize vehicle's inclination angle}%
\includegraphics[width=8.5 cm, height=5.5cm]{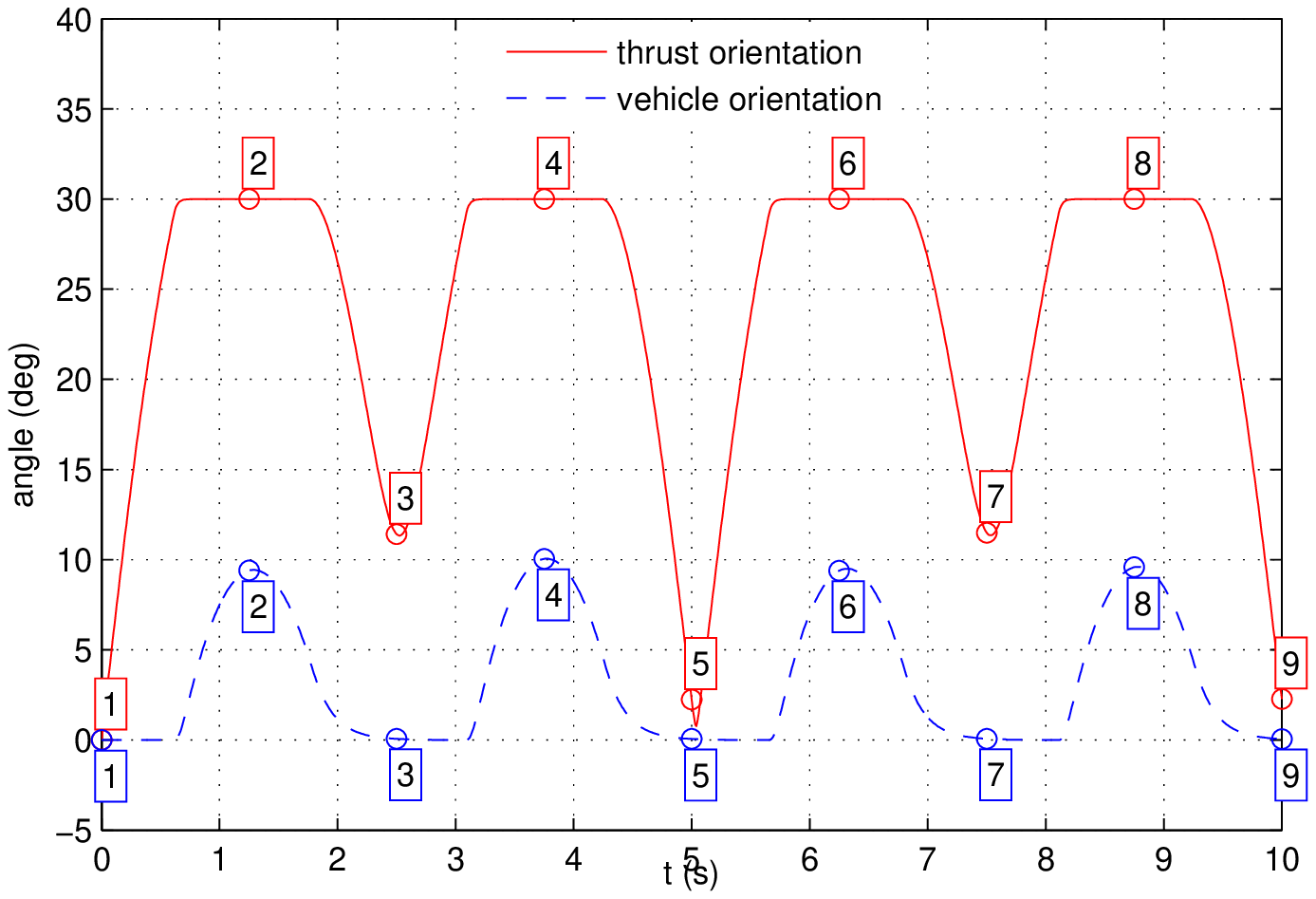}
\vspace*{-0.3cm}
\caption{Thrust and vehicle inclinations vs. time (Simulation 2).} \label{fig:simu2b} \vspace*{0.2cm}
\psfrag{t (s)}{\scriptsize $t\, (s)$}%
\psfrag{position error (m)}{\scriptsize $Position$ $errors$ $\tilde{x}$ $(m)$}%
\psfrag{error x}{\scriptsize $\tilde{x}_1(m)$}%
\psfrag{error y}{\scriptsize $\tilde{x}_2(m)$}%
\psfrag{error z}{\scriptsize $\tilde{x}_3(m)$}%
\psfrag{vehicle orientation}{\scriptsize vehicle's inclination angle}%
\includegraphics[width=8.5 cm, height=5.5cm]{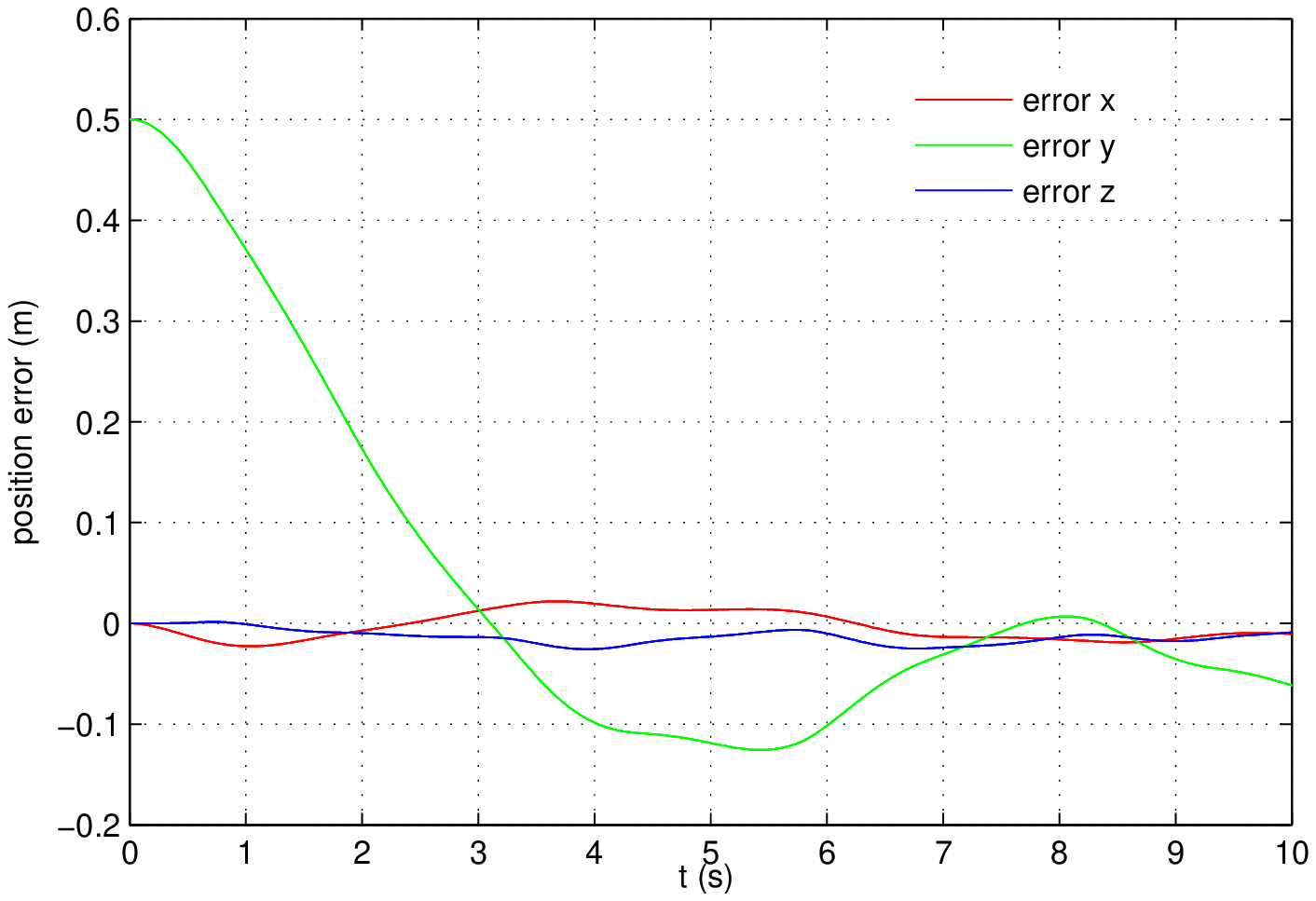}
\vspace*{-0.3cm}
\caption{Position tracking errors vs. time (Simulation 2).} \label{fig:simu2c}
\end{figure}

\begin{figure}
\psfrag{t (s)}{\scriptsize $t\, (s)$}%
\psfrag{Thrust (N)}{\scriptsize $Thrust$ $T(N)$}%
\includegraphics[width=8.5 cm, height=5.5cm]{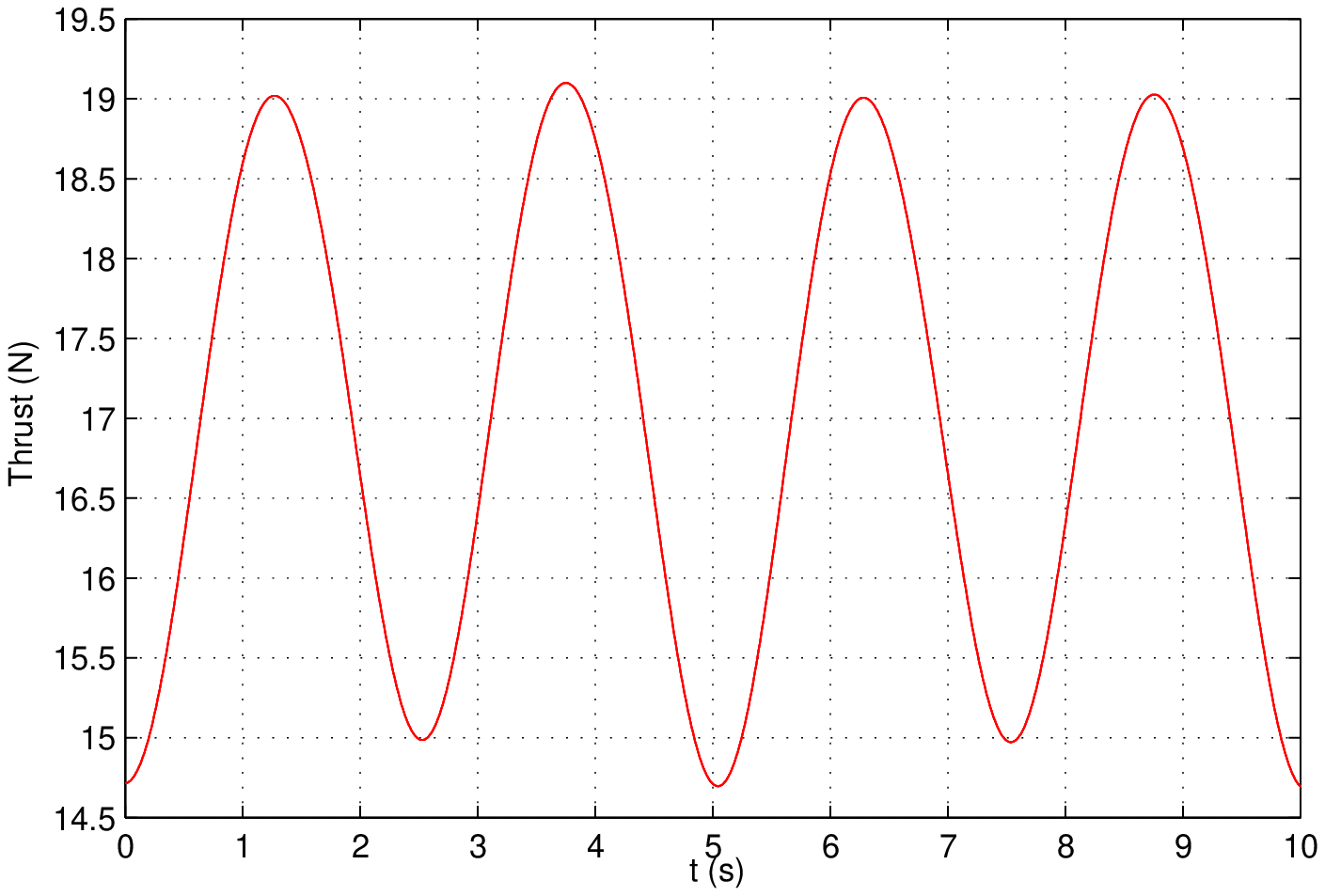}
\vspace*{-0.3cm}
\caption{Thrust intensity vs. time (Simulation 2).} \label{fig:simu2d} \vspace*{0.2cm}
\psfrag{t (s)}{\scriptsize $t\, (s)$}%
\psfrag{Torque (Nm)}{\scriptsize $Torque$  $\Gamma(Nm)$}%
\psfrag{gamma x}{\tiny $\Gamma_1 (Nm)$ }%
\psfrag{gamma y}{\tiny $\Gamma_2 (Nm)$}%
\psfrag{gamma z}{\tiny $\Gamma_3 (Nm)$}%
\includegraphics[width=8.5 cm, height=5.5cm]{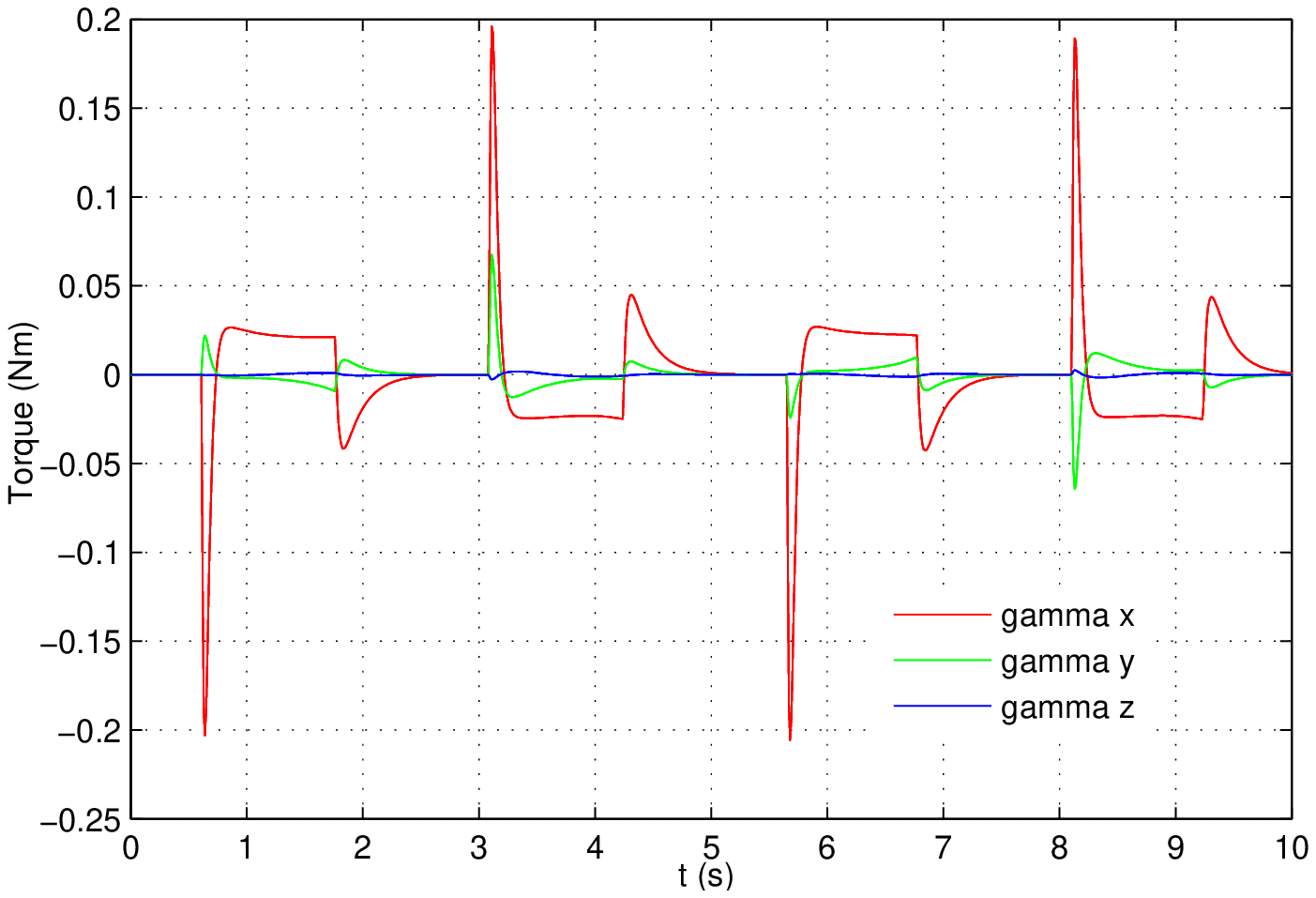}
\vspace*{-0.3cm}
\caption{Torque inputs vs. time (Simulation 2).} \label{fig:simu2e}
\end{figure}

\begin{itemize}
\item {\em Simulation 1:} The time period for a complete run is $15 s$. Beside the projection on the horizontal plane of the path followed by the vehicle's CoM (Fig. \ref{fig:simu1a}), the time evolution of the vehicle's and thrust-direction's inclinations, of the position tracking errors, and of the thrust magnitude, is shown on Figs. \ref{fig:simu1b}--\ref{fig:simu1d}. In Figs. \ref{fig:simu1a} and \ref{fig:simu1b}, the nine highlighted points correspond to time-instants when the reference trajectory involves large acceleration variations. From Fig. \ref{fig:simu1a}, one can observe that the vehicle catches up and subsequently closely follows the reference trajectory. Despite a rather aggressive reference trajectory (with an average longitudinal velocity of about $4 m/s$ and accelerations sometimes exceeding $3 m/s^2$), the vehicle's base remains always horizontal (see Fig. \ref{fig:simu1b}). One can also observe from Fig. \ref{fig:simu1b} that the thrust-direction tilt angle never reaches its maximal value (equal to $\pi/6 rad$). Both objectives are thus realized (almost) perfectly in this case.

\item {\em Simulation 2:} This simulation is devised to illustrate the effects of thrust tilting saturation and the corresponding control monitoring. The time-period for a complete run is reduced to $10s$ and results are shown in Figs. \ref{fig:simu2a}--\ref{fig:simu2e}. One can now observe from Fig. \ref{fig:simu2b} that the vehicle's inclination periodically departs from zero when the thrust-direction tilt angle attains its maximal value. While the position tracking primary objective remains perfectly realized, the secondary objective is not realized in this case. However, the body's inclination returns to the desired zero value as soon as the thrust tilt angle that is required to achieve the second objective re-enters the domain of allowed tilt angles, a behaviour that we find satisfactory.
\end{itemize}


\section{Conclusion} \label{conclusion}
Nonlinear control of VTOL vehicles endowed with thrust tilting capability has been addressed and a generic solution to the thrust-tilting problem has been devised.
The proposed solution potentially applies to a large panel of aerial vehicles with extended flight envelopes. It involves a primary objective consisting in the asymptotic stabilization of either a reference velocity or a reference position trajectory, and a secondary objective consisting in the asymptotic stabilization of either a reference direction for one of the body-base vectors or of a complete reference attitude for the body-fixed frame. A major original outcome of the present study is the definition of a control solution that takes thrust-tilting limitations into account explicitly. We view it also as new contribution to the ongoing development of the unified approach to the control of aerial vehicles. The validity of the proposed thrust tilting control strategy has so far been validated only in simulation and we are now interested in applications and extensions carried out on physical devices.

\end{document}